\documentclass[journal]{new-aiaa}
\usepackage[utf8]{inputenc}

\usepackage{graphicx}
\usepackage{amsmath}
\usepackage[version=4]{mhchem}
\usepackage{siunitx}
\usepackage{longtable,tabularx}
\usepackage{subcaption}
\usepackage{todonotes}
\usepackage{bm}
\usepackage{multirow}
\usepackage{xcolor}
\usepackage{comment}

\newcommand{\R}{\mathbb{R}}
\title{Multidisciplinary Design Optimization Approach to \\ Integrated Space Mission Planning and Spacecraft Design\footnote{This paper is a revised version of Paper 2021-4069 at the AIAA ASCEND 2021, Las Vegas, Nevada \& Virtual, November 15-17, 2021.}}

\author{Masafumi Isaji\footnote{Ph.D. Student, Daniel Guggenheim School of Aerospace Engineering, Atlanta, GA, AIAA Student Member.}, Yuji Takubo\footnote{Undergraduate Student, Daniel Guggenheim School of Aerospace Engineering, Atlanta, GA, AIAA Student
Member.}, and Koki Ho\footnote{Assistant Professor, Daniel Guggenheim School of Aerospace Engineering, Atlanta, GA, AIAA Senior Member.}}
\affil{Georgia Institute of Technology, Atlanta, GA, 30332, USA}

\begin{document}

\maketitle

\begin{abstract}
Space mission planning and spacecraft design are tightly coupled and need to be considered together for optimal performance; however, this integrated optimization problem results in a large-scale Mixed-Integer Nonlinear Programming (MINLP) problem, which is challenging to solve. In response to this challenge, this paper proposes a new solution approach to this problem based on decomposition-based optimization via augmented Lagrangian coordination. The proposed approach leverages the unique structure of the problem that enables its decomposition into a set of coupled subproblems of different types: a Mixed-Integer Quadratic Programming (MIQP) subproblem for mission planning, and one or more Nonlinear Programming (NLP) subproblem(s) for spacecraft design. Since specialized MIQP or NLP solvers can be applied to each subproblem, the proposed approach can efficiently solve the otherwise intractable integrated MINLP problem. An automatic and effective method to find an initial solution for this iterative approach is also proposed so that the optimization can be performed without a user-defined initial guess.
The demonstration case study shows that, compared to the state-of-the-art method, the proposed formulation converges substantially faster and the converged solution is also at least the same or better given the same computational time limit.

\end{abstract}

\section{Nomenclature}

{\renewcommand\arraystretch{1.0}
\noindent\begin{longtable*}{lll}
$\mathcal{A}$ &\quad=\quad& Set of arcs \\
$\boldsymbol{a}_{vijt}$ &\quad=\quad& Cost coefficient vector of commodity\\
${a'}_{vijt}$ &\quad=\quad& Cost coefficient of spacecraft\\
$\boldsymbol{d}_{it}$ &\quad=\quad& Demand vector\\
$\boldsymbol{e}_v$ &\quad=\quad& Spacecraft design variable vector \\
$\mathcal{F}(-)$ &\quad=\quad& Spacecraft sizing function\\
$f$ &\quad=\quad& Objective function (subproblem) \\
$\boldsymbol{g}$ &\quad=\quad& Inequality constraint \\
$\boldsymbol{h}$ &\quad=\quad& Equality constraint \\
$H_{vij}$ &\quad=\quad& Concurrency matrix \\
$\mathcal{J}$ &\quad=\quad& Objective function \\
$L$ &\quad=\quad& Number of subsystems in the dry mass\\
$M$ &\quad=\quad& Number of subproblems in a quasi-separable MDO problem\\
$m$ &\quad=\quad& Mass of spacecraft subsystems \\
$m_d$ &\quad=\quad& Spacecraft dry mass \\
$m_f$ &\quad=\quad& Spacecraft propellant capacity \\
$m_p$ &\quad=\quad& Spacecraft payload capacity \\
$N$ &\quad=\quad& Number of types of spacecraft \\
$\mathcal{N}$ &\quad=\quad& Set of nodes \\ 
$n$ &\quad=\quad& Dimension of variables \\
$Q_{vijt}$ &\quad=\quad& Commodity transformation matrix \\
$\mathcal{T}$ &\quad=\quad& Set of time steps \\
$t_{mis}$ &\quad=\quad& Mission length \\
$\Delta t_{ij}$ &\quad=\quad& Time of Flight (ToF)\\
$u_{vijt}$ &\quad=\quad& Spacecraft flow variable \\
$\mathcal{V}$ &\quad=\quad& Set of spacecraft \\
$W_{ij}$ &\quad=\quad& Launch time window \\ 
$\boldsymbol{x}_{vijt}$ &\quad=\quad& Commodity flow variable\\
$\boldsymbol{y}$ &\quad=\quad& Shared variables \\
$\boldsymbol{z}$ &\quad=\quad& Local variables \\
$\zeta$ &\quad=\quad& Propellant type \\
$\phi$ &\quad=\quad& Penalty function \\ \\

\emph{Subscipt \& Superscript}&&\\ 
$i$ &\quad=\quad& Node index (departure) \\
$j$ &\quad=\quad& Node index (arrival) \\
$k$ &\quad=\quad& Subproblem index  \\
$l$ &\quad=\quad& Subsystem index  \\
$q$ &\quad=\quad& Iteration count\\
$t$ &\quad=\quad& Time index \\
$v$ &\quad=\quad& Vehicle index \\
\end{longtable*}}
\section{Introduction}

\lettrine{A}{s} we pursue sustainable presence in space, a framework to optimize large-scale, long-term space missions efficiently is imperative. A number of studies on space logistics that incorporate the transportation network in large-scale space mission design have been developed, including SpaceNet \cite{Shull2007MS}, the interplanetary logistics model \cite{taylor2007logistics}, and the extensive literature on space logistics optimization frameworks based on the generalized multicommodity network flow \cite{ishimatsu2016gmcnf, ho2014time-expanded, ho2016FlexiblePath}. Utilizing the linear nature of such space logistics or transportation network optimization problems, researchers have developed frameworks that can efficiently optimize the mission design as Mixed-Integer Linear Programming (MILP) problems \cite{chen2018MILP, chen2018regular, chen2019isruMars, takubo2021HRL}. 
However, due to the nonlinear nature of spacecraft design, a naive integration of spacecraft design into space mission/campaign planning (a transportation scheduling or resource distribution) would result in a large-scale Mixed-Integer Nonlinear Programming (MINLP) problem. Even though this integration has been shown to exhibit substantial solution improvement compared to optimizing them individually \cite{taylor2007phd}, solving the MINLP is oftentimes computationally prohibitive. 
Since the concurrent optimization of space mission planning and spacecraft design is highly desired in practice, each community took different approaches to bridge these two domains.


In the space logistics community, spacecraft design has been considered as a high-level nonlinear sizing model and has been integrated into mission planning either by separating the nonlinear part from the mission planning optimization or by piecewise linearization of the spacecraft model. Taylor \cite{taylor2007phd} developed a parametric spacecraft sizing model which determines the spacecraft dry mass from its payload capacity and propellant capacity. Based on this model, Simulated Annealing (SA) or a similar metaheuristic optimization algorithm optimizes the spacecraft design variables, while the linear programming (LP) or MILP solver evaluates the constraints and determines transportation flow variables. In this way, the LP or MILP solver is embedded into SA, and thus it was called the embedded optimization methodology.
Using the same spacecraft sizing model, Chen and Ho \cite{chen2018MILP} employed the piecewise linear (PWL) approximation of the nonlinear model to approximate the entire MINLP problem as a MILP problem that can be solved efficiently. However, this approach is an approximation model, and the resulting solution is not guaranteed to be feasible nor optimal in the original nonlinear problem. For a different yet related problem, satellite component selection and operation have been co-optimized as a mixed-integer programming problem \cite{norheim2021, norheim2020}; however, the existing techniques are based on the specific types of nonlinear relationships in the problem (e.g., power laws) and are not generally applicable to the integrated logistics problem.

On the other hand, aerospace vehicle design has been tackled by the Multidisciplinary Design Optimization (MDO) community. Despite various optimization and sizing methods that can deal with the high-dimensional nonlinear design of aircraft or spacecraft \cite{sobieszczanski1997MDOsurvey}, few studies integrated the mission-level analysis or optimization. One of the few studies that tackled the integrated mission planning and spacecraft design is Ref. \cite{beauregard2021lunarMDO} by Beauregard et al., which proposed an MDO architecture for a lunar lander design with a lunar mission sequence architecture analysis. This architecture connects the mission planning and spacecraft design problem using a sequential procedure without a feedback structure (i.e., the mission architecture is first chosen and fixed, then the lunar lander MDO is performed); therefore, the mission and spacecraft are not simultaneously optimized, and spacecraft design is neglected when selecting the mission architecture. In addition, the candidates of the mission architectures are given \emph{a priori} and discrete (combinatory). These two factors limit the design space and make this approach not suitable for the integrated space mission design. 

This paper proposes an efficient decomposition-based optimization scheme for integrated space mission planning and spacecraft design. The key idea is to decompose the integrated MINLP problem into multiple coupled subproblems of different types: the Mixed-Integer Quadratic Programming (MIQP) subproblem for space mission planning, and the  Nonlinear Programming (NLP) subproblem(s) for spacecraft design. Since specialized efficient MIQP or NLP optimizers (e.g., Gurobi \cite{gurobi} for MIQP; IPOPT \cite{IPOPT} for NLP) can be utilized to solve each subproblem, the proposed method can solve the otherwise intractable integrated MINLP problem efficiently.
The iterative coordination between each subproblem can be achieved using an MDO approach \cite{martins2013MDOsurvey,sobieszczanski1997MDOsurvey}. Specifically, the Augmented Lagrangian Coordination (ALC) approach \cite{tosserams2007ALC} with the Analytical Target Cascading (ATC) structure \cite{ATC,ATCextended} is chosen for the proposed method. This architecture fits our problem well because (1) it allows us to decompose the original complex problem into the subproblems with different and simpler types (MIQP or NLP), each of which can be efficiently solvable with specialized solvers; (2) it has a robust convergence property; and (3) it allows the complex hierarchical structure for the spacecraft design subproblem(s) and can be easily parallelizable (and thus scalable) if needed. Since the nonlinear optimization solvers generally require a good initial guess, we further develop an automated initial guess generation method based on PWL approximation to the MINLP problem so that no user-defined initial guess is needed for the optimization. The demonstration case study shows that, compared to the state-of-the-art method, the proposed new formulation converges substantially faster and the converged solution is also at least the same or better given the same computational time limit. 

The remainder of this paper proceeds as follows. In Section \ref{Problem Definition}, the problem definition of the integrated space mission planning and spacecraft design as an all-in-one optimization problem formulation is described.  Section \ref{OurALC} illustrates the solution procedure for the proposed problem based on the decomposition-based method. Section \ref{CaseStudy} introduces a case study of human lunar exploration missions and compares the computational efficiency of the proposed method against the state-of-the-art method. Finally, Section \ref{conclusion} states the conclusion.

\section{Problem Definition: Integrated Space Mission Planning and Spacecraft Design} \label{Problem Definition}
The goal of this research is to optimize the transportation scheduling (referred to as space mission planning) and vehicle design (referred to as spacecraft design) for a long-term space campaign that can potentially comprise multiple missions. This section introduces the formulation for this integrated space mission planning and spacecraft design problem (referred to as the all-in-one formulation). The idea behind this formulation is to consider space mission planning as a transportation network optimization problem for which the design of vehicles is also part of the decision variables. In the network, the nodes correspond to the orbital or surface locations, and the arcs correspond to the trajectories connecting the nodes. The decision variables include both the commodities that flow over the network and the design parameters for the vehicles that carry these commodities. The optimization formulation is listed as follows, and the list of variables and parameters is included in Table~\ref{tab_SLvar}. 

\begin{equation}
\label{SL_obj}
\min \quad \mathcal{J}=\sum_{t\in \mathcal{T}} \sum_{(v,i,j)\in \mathcal{A}} (\boldsymbol{a}_{vijt}^{T} \boldsymbol{x}_{vijt} + {a'}_{vijt}^{T} m_{d_{v}} u_{vijt})
\end{equation}

subject to
\begin{equation}
\label{SL_constr1}
\sum_{(v,j):(v,i,j)\in \mathcal{A}} \left[\begin{array}{c}
\boldsymbol{x}_{vijt} \\
m_{d_v}u_{vijt}
\end{array}\right] - \sum_{(v,j):(v,i,j)\in \mathcal{A}} Q_{vjit} \left[\begin{array}{c}
\boldsymbol{x}_{vji(t-\Delta t_{ji})} \\
m_{d_{v}} u_{vji(t-\Delta t_{ji})}
\end{array}\right]
\leq \boldsymbol{d}_{it} \quad \forall t\in \mathcal{T}\quad \forall i\in \mathcal{N} 
\end{equation}

\begin{equation}
\label{SL_constr2}
    {H}_{vij} \boldsymbol{x}_{vijt} \leq \boldsymbol{e}_v u_{vijt}  \quad \forall t\in \mathcal{T} \quad \forall (v,i,j)\in \mathcal{A}
\end{equation}

\begin{equation}
\label{SL_constr3}
\left\{\begin{array}{ll}
\boldsymbol{x}_{vijt} \geq \boldsymbol{0}_{p \times 1}  &\text { if } t \in W_{i j} \\
\boldsymbol{x}_{vijt} =\boldsymbol{0}_{p \times 1} \quad &\text { otherwise }
\end{array} \quad \forall(v, i, j) \in \mathcal{A}\quad \forall t\in \mathcal{T}\right.
\end{equation}

\begin{equation}
\label{SL_vehicle_sizing}
    m_{d_v} = \mathcal{F}(\boldsymbol{e}_v, \zeta_v)  \quad \forall v\in \mathcal{V} 
\end{equation}

\begin{equation}
\label{xdef}
\begin{aligned}
\boldsymbol{x}_{v i j t}=\left[\begin{array}{c}
x_{1} \\
x_{2} \\
\vdots \\
x_{p}
\end{array}\right]_{v i j t}, \quad \begin{aligned}
&x_{n} \in \mathbb{R}_{\geq 0} \quad \forall n \in \mathcal{C}_{c} \\
&x_{n} \in \mathbb{Z}_{\geq 0} \quad \forall n \in \mathcal{C}_{d}
\end{aligned} \quad \forall(v, i, j, t) \in \mathcal{A} \\
\end{aligned}
\end{equation}

\begin{equation}
\label{udef}
{u}_{v i j t} \in \mathbb{Z}_{\geq 0}  \quad \forall(v, i, j, t) \in \mathcal{A} 
\end{equation}

\begin{equation}
\label{evdef}
\boldsymbol{e}_v =  \left[\begin{array}{c}
m_{p} \\
m_{f} 
\end{array}\right]_{v}, \quad m_{p_v}, m_{f_v}, m_{d_v} \in \mathbb{R}_{\geq 0}, \quad \forall v \in \mathcal{V} 
\end{equation}

\begin{table}[h]
\caption{\label{tab_SLvar} Variables and parameters used in the space transportation scheduling problem}
\centering
\begin{tabular}{ll}
\hline\hline Name & 
\begin{tabular}{p{0.7\textwidth}}Description\end{tabular} \\\hline

{\emph{Variables} }&\\

\hline
$\boldsymbol{x}_{vijt}$ & 
\begin{tabular}{p{0.7\textwidth}}
Commodity flow variable, or the quantity of the commodity delivered from node $i$ to $j$ at time $t$ by spacecraft $v$. $\boldsymbol{x}_{vijt} \geq 0$. Each component of this variable can contain either continuous variables ($\mathcal{C}_{c}$) or discrete variables ($\mathcal{C}_{d}$). This vector will be $\mathbb{R}^p$ if the total commodity variation is $p$. 
\end{tabular} \\

$u_{vijt}$ & 
\begin{tabular}{p{0.7\textwidth}}
Spacecraft flow variable, which indicates the number of spacecraft type $v$ moving from node $i$ to $j$ at time $t$ ($\in \mathbb{R}$). This variable is integer scalar. 
\end{tabular} \\

$\boldsymbol{e}_v$ & 
\begin{tabular}{p{0.7\textwidth}}
Spacecraft design variables and parameters. In this problem, it includes payload capacity $m_p$ and propellant capacity $m_f$ ($\in \mathbb{R}^2$).
\end{tabular}\\

$m_{d_v}$ & 
\begin{tabular}{p{0.7\textwidth}}
Dry mass of spacecraft $v$ ($\in \mathbb{R}$).
\end{tabular}\\

\hline
{\emph{Parameters}} &\\

\hline
$\boldsymbol{a}_{vijt}$ & 
\begin{tabular}{p{0.7\textwidth}}
Cost coefficient of commodity ($\in \mathbb{R}^p$).
\end{tabular}\\

${a'}_{vijt}$ & 
\begin{tabular}{p{0.7\textwidth}}
Cost coefficient of spacecraft ($\in \mathbb{R}$). 
\end{tabular}\\

$\boldsymbol{d}_{it}$ & 
\begin{tabular}{p{0.7\textwidth}}
Demands/supplies of different commodities and spacecraft at node $i$ at time $t$ ($\in \mathbb{R}^{p+1}$).
\end{tabular}\\

$Q_{vijt}$ & 
\begin{tabular}{p{0.7\textwidth}}
Transformation matrix ($\in\R^{(p+1)\times(p+1)}$).
\end{tabular}\\

$H_{vij}$ & 
\begin{tabular}{p{0.7\textwidth}}
Concurrency constraint matrix ($\in \R^{2\times p}$).
\end{tabular}\\

$W_{ij}$ & 
\begin{tabular}{p{0.7\textwidth}}
Launch window vector, which indicates the available launch window of spacecraft.
\end{tabular}\\

$\mathcal{F}(-)$ & 
\begin{tabular}{p{0.7\textwidth}}
Spacecraft sizing function. This illustrates the nonlinear relationship of the spacecraft design variables and design parameters. 
\end{tabular}\\

$\Delta t_{ij}$ & 
\begin{tabular}{p{0.7\textwidth}}
Time of Flight (ToF) from node $i$ to $j$.
\end{tabular}\\

$\zeta_v$ & 
\begin{tabular}{p{0.7\textwidth}}
Propellant type for each spacecraft (predetermined). 
\end{tabular}\\

\hline
\emph{Sets} &\\
\hline
$\mathcal{A(V,N,N,T)}$ & 
\begin{tabular}{p{0.7\textwidth}}
Set of arcs realized by spacecraft.
\end{tabular}\\

$\mathcal{N}$ & 
\begin{tabular}{p{0.7\textwidth}}
Set of nodes.
\end{tabular}\\

$\mathcal{T}$ & 
\begin{tabular}{p{0.7\textwidth}}
Set of time steps.
\end{tabular}\\

$\mathcal{V}$ & 
\begin{tabular}{p{0.7\textwidth}}
Set of spacecraft (vehicles).
\end{tabular}\\

\hline\hline
\end{tabular}
\end{table}

Equation~\eqref{SL_obj} indicates the objective function, which can be the lifecycle cost or launch mass, depending on the application context. In this research, we set the coefficients $\boldsymbol{a}_{vijt}$ and ${a}_{vijt}$ so that the objective function corresponds to the sum of initial mass at low-earth orbit (IMLEO), a metric commonly used in the space logistics literature as can be seen in Ref. \cite{taylor2007phd,ishimatsu2016gmcnf,ho2014time-expanded,chen2018MILP}.

Equations~\eqref{SL_constr1}-\eqref{SL_constr3} are the constraints for space mission planning.
First, Eq.~\eqref{SL_constr1} is the mass balance constraint that guarantees that the inflow (supply) of the commodity is larger than the sum of the outflow and demand. $Q_{vijt}$ is the transformation matrix, which indicates the transformation of the commodity during the spaceflight; for example, the relationship of impulsive propellant consumption can be illustrated using this constraint. 
Next, Eq.~\eqref{SL_constr2} is the concurrency constraint. This indicates that the commodity loaded on each spacecraft is constrained by the dimension of the spacecraft. Specifically, in this paper, the payload and propellant flow is limited: the amount of propellant is lower than the propellant capacity of the spacecraft, and the sum of other payloads is lower than the payload capacity.
Finally, Eq.~\eqref{SL_constr3} is the time window constraints. The commodity flow is allowed only if the time $t$ belongs to the launch window vector $W_{ij}$, and for the remaining time steps, the commodity flow is conserved to be zero. 

Equation~\eqref{SL_vehicle_sizing} indicates an abstract representation of the spacecraft design constraints, which describes the constraints between the properties of the vehicle. It can take a wide range of complexity, including an explicit or implicit relationship of the subsystems or design parameters of the spacecraft; when the spacecraft requires multiple disciplines or multiple subsystems, an MDO problem can be embedded in this constraint. 

Along with Table~\ref{tab_SLvar}, Eqs. \eqref{xdef}, \eqref{udef}, and \eqref{evdef} show the definitions and domains of commodity flow variables, spacecraft flow variables, and spacecraft design variables, respectively.

This integrated mission planning and spacecraft design problem results in a constrained MINLP problem, one of the most challenging optimization problem types to solve. Namely, this problem contains both discrete and continuous variables as well as both linear and nonlinear constraints. In particular, the discrete variables represent the definition of the commodity flow and the number of spacecraft on the mission planning side of the problem. The nonlinearity appears in two ways: (1) the spacecraft design relationship in Eq.~\eqref{SL_vehicle_sizing}; (2) the quadratic terms in the mass balance constraint (Eq.~\ref{SL_constr1}) and concurrency constraint (Eq.~\ref{SL_constr2}) for mission planning (note: both $\boldsymbol{e}_v$ and $u_{vijt}$ are variables). Fortunately, this second nonlinearity can be converted into an equivalent linear relationship through the big-M method, as explained in Ref.~\cite{chen2018MILP}, so that the nonlinearity only exists on the spacecraft design side of the problem.
Therefore, as a result, the problem contains two coupled problems: one for space mission planning which is linear with integer variables, and the other for spacecraft design which is nonlinear with continuous variables.
Our approach leverages this unique structure of the problem and proposes a new approach to solve this problem efficiently.

\section{Proposed Approach: Decomposition-Based Optimization with Augmented Lagrangian Coordination}
\label{OurALC}
\begin{figure}[h]	
	\centering
	\begin{subfigure}[t]{1.8in}
		\centering
		\includegraphics[scale=0.5]{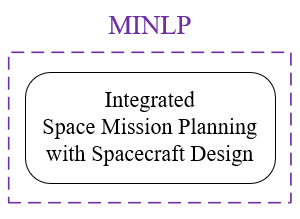}
		\subcaption{All-in-one formulation}\label{fig:1a}		
	\end{subfigure}
	\begin{subfigure}[t]{4.5in}
		\centering
		\includegraphics[scale=0.36]{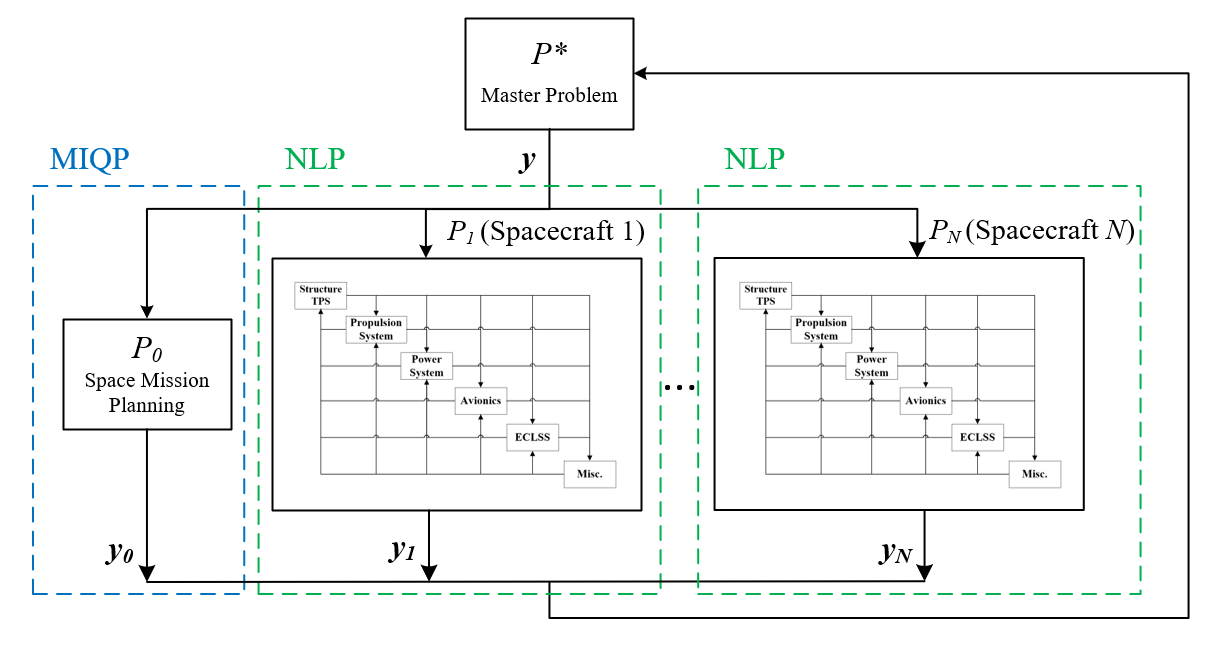}
		\subcaption{Proposed decomposition-based formulation based on \cite{tosserams2007ALC}}\label{fig:1b}
	\end{subfigure}
	\caption{Solution strategy for integrated space mission planning and spacecraft design.}\label{fig:1}
\end{figure}

Decomposition-based optimization is often used to decompose an MDO problem in terms of disciplines or subsystems. Exploiting the unique feature of the integrated space mission planning and spacecraft design problem, we apply this approach to decompose the large-scale MINLP problem (Fig.~\ref{fig:1a}) into coupled MIQP and NLP subproblems (Fig.~\ref{fig:1b}), each of which is significantly easier to solve with specialized solvers compared to the original MINLP problem. The space mission planning subproblem can be solved using a MIQP solver, and the spacecraft design subproblem can be solved using an NLP solver without any integer variables. The coupled subproblems are solved iteratively using the ALC-based coordination until convergence is reached. To enable the optimization without a user-defined initial guess, an automated and effective initial solution generation approach is also proposed.

\subsection{Derivation of Decomposed Problems with Augmented Lagrangian Coordination}
\label{OriginalALC}
We first start with deriving the formulations of the decomposed problems with ALC. ALC tackles complex MDO optimization problems that are quasi-separable and thus can be decomposed into a set of coupled subproblems. ALC is attractive because of (1) its ability to break down our MINLP problem into MIQP and NLP problems; (2) its robust convergence property; and (3) its flexibility with the hierarchical structure of the problems. For an extensive discussion on ALC, refer to Ref.~\cite{tosserams2007ALC}.

The formulation for the quasi-separable MDO problem with $M$ subproblems is given as follows:

\begin{equation}
\begin{array}{rl}
\underset{\boldsymbol{y}, \boldsymbol{z}_{0}, \ldots, \boldsymbol{z}_{M-1}} {\text{ min }} \quad
\displaystyle{\sum_{k=0}^{M-1} f_{k}\left(\boldsymbol{y}, \boldsymbol{z}_{k}\right)} \\

\text { subject to } 
\quad \boldsymbol{g}_{k}\left(\boldsymbol{y}, \boldsymbol{z}_{k}\right) \leq \boldsymbol{0} & k=0, \ldots, M-1 \\
\quad \boldsymbol{h}_{k}\left(\boldsymbol{y}, \boldsymbol{z}_{k}\right)=\boldsymbol{0} & k=0, \ldots, M-1
\end{array}
\end{equation}

\noindent where $\boldsymbol{y} \in \mathbb{R}^{n^{y}}$ indicates the shared variables, $\boldsymbol{z}_k \in \mathbb{R}^{n_{k}^{z}}$ indicates the local variables for subproblem $k$. The shared variables $\boldsymbol{y}$ can be common variables over multiple subproblems. $f_k:\mathbb{R}^{n_{k}} \mapsto \mathbb{R}$ indicates the local objective function, $\boldsymbol{g}_k$ and $\boldsymbol{h}_k$ indicate the equality and inequality constraints for each subproblem. The dimension of the total design variable  $\boldsymbol{s}=\left[\boldsymbol{y}^{T}, \boldsymbol{z}_{0}^{T}, \ldots, \boldsymbol{z}_{M-1}^{T}\right]^{T} \in \mathbb{R}^{n}$ is $n = n^{y} + \sum_{k=0}^{M-1} n_{k}^{z}$. The dimension of the design variable for subproblem $k$ is $n_k = n^{y} + n_{k}^{z}$.  

The decomposition-based approach for this problem follows the following steps. First, we introduce the auxiliary variables and consistency constraints so that the local constraints, $\boldsymbol{g}_k$ and $\boldsymbol{h}_k$, are only dependent on the auxiliary variables $\boldsymbol{y}_k$ and independent of the shared variables $\boldsymbol{y}$. 

\begin{equation}
\begin{alignedat}{2}
\min_{\boldsymbol{y}, \boldsymbol{y}_{0}, \boldsymbol{z}_{0}, \ldots, \boldsymbol{y}_{M-1}, \boldsymbol{z}_{M-1}} \quad &\sum_{k=0}^{M-1} f_{k}\left(\boldsymbol{y}_{k}, \boldsymbol{z}_{k}\right) & \\
\text { subject to } 
\quad &\boldsymbol{g}_{k}\left(\boldsymbol{y}_{k}, \boldsymbol{z}_{k}\right) \leq \boldsymbol{0} &\quad k=0, \ldots, M-1 \\
&\boldsymbol{h}_{k}\left(\boldsymbol{y}_{k}, \boldsymbol{z}_{k}\right)=\boldsymbol{0} &\quad k=0, \ldots, M-1 \\
&\boldsymbol{c}_{k}\left(\boldsymbol{y}, \boldsymbol{y}_{k}\right)=\boldsymbol{0} &\quad k=0, \ldots, M-1 
\end{alignedat}
\end{equation}

\noindent With the consistency constraints $\boldsymbol{c}_{k}$, which ensures that the auxiliary variables $\boldsymbol{y}_{k}$ are the same as the shared variables $\boldsymbol{y}$, the shared variables are separated from the local variables while representing the same problem as the original one. Next, the relaxation of the consistency constraints is introduced with the local Lagrangian penalty function:

\begin{equation}
\begin{alignedat}{1}
\min_{\boldsymbol{y}, \boldsymbol{y}_{0}, \boldsymbol{z}_{0}, \ldots, \boldsymbol{y}_{M-1}, \boldsymbol{z}_{M-1}} \quad &\sum_{k=0}^{M-1} f_{k}\left(\boldsymbol{y}_{k}, \boldsymbol{z}_{k}\right) +\sum_{k=0}^{M-1} \phi_{k}\left(\boldsymbol{c}_{k}\left(\boldsymbol{y}, \boldsymbol{y}_{k}\right)\right) \\
\text { subject to } 
\quad &\boldsymbol{g}_{k}\left(\boldsymbol{y}_{k}, \boldsymbol{z}_{k}\right) \leq \boldsymbol{0} \quad k=0, \ldots, M-1 \\
\quad &\boldsymbol{h}_{k}\left(\boldsymbol{y}_{k}, \boldsymbol{z}_{k}\right)=\boldsymbol{0} \quad k=0, \ldots, M-1
\end{alignedat}
\end{equation}

\noindent The augmented Lagrangian penalty function for subproblem $k$, $\phi_k$, is defined as follows. 

\begin{equation}
\phi_{k}\left(\boldsymbol{c}_{k}\left(\boldsymbol{y}, \boldsymbol{y}_{k}\right)\right)=\boldsymbol{v}_{k}^{T}\left(\boldsymbol{y}-\boldsymbol{y}_{k}\right)+\left\|\boldsymbol{w}_{k} \circ\left(\boldsymbol{y}-\boldsymbol{y}_{k}\right)\right\|_{2}^{2}
\end{equation}

\noindent where $\boldsymbol{v_k}$ is the vector of Lagrange multiplier estimates, and $\boldsymbol{w_k}$ is the vector of penalty weights. Here, $\circ$ represents the element-wise product of matrices or vectors, also known as the Hadamard product. By moving the consistency constraints into the local objective functions, the local subproblems can be completely separated.
The bi-level decomposition-based problem is now formulated by establishing the master problem above the subproblems. The master problem minimizes the penalty function and updates the shared variables $\boldsymbol{y}$. Note that even though the bi-level formulation is employed here, the ALC has the capability to handle multi-level hierarchical formulation as well \cite{ATCextended}.

(1) Master Problem
\begin{equation}
\min _ {\boldsymbol{y}} \quad \sum_{k=0}^{M-1} \phi_{k}\left(\boldsymbol{c}_{k}\left(\boldsymbol{y}, \boldsymbol{y}_{k}\right)\right)
\end{equation}

(2) Subproblem $k$
\begin{equation}
\begin{aligned}
\underset{\boldsymbol{y}_k, \boldsymbol{z}_{k}} {\text{ min }} \quad &f_{k}\left(\boldsymbol{y}_{k}, \boldsymbol{z}_{k}\right) +\phi_{k}\left(\boldsymbol{c}_{k}\left(\boldsymbol{y}, \boldsymbol{y}_{k}\right)\right) \\
\text { subject to } 
\quad&\boldsymbol{g}_{k}\left(\boldsymbol{y}_{k}, \boldsymbol{z}_{k}\right) \leq \boldsymbol{0} \\
&\boldsymbol{h}_{k}\left(\boldsymbol{y}_{k}, \boldsymbol{z}_{k}\right) =\boldsymbol{0}
\end{aligned}
\end{equation}

Adopting the above approach in our problem of the integrated space mission planning and spacecraft design with $N$ vehicle types, Fig.~\ref{fig:1b} represents the decomposition-based optimization architecture. We have one space mission planning subproblem (Subproblem $0$) and multiple spacecraft design subproblems (Subproblems $1,\ldots,N$), where $N$ is the number of spacecraft types; thus, we have $N+1$ subproblems in total (i.e., $M=N+1$). The shared variables among them include the vehicle design parameters $\boldsymbol{y} = [\boldsymbol{y}_1^T, \ldots, \boldsymbol{y}_N^T]^T$, where $\boldsymbol{y}_v = [m_{p_v}, m_{f_v}, m_{d_v}]^T$, where $v$ is the vehicle index such that $v = 1,\ldots, N$. Note that, in the considered problem, the subproblem index $k$ and vehicle index $v$ match (i.e., subproblem $k$ handles the design problem of spacecraft $v$). For each spacecraft $v$, $m_{p_v}$, $m_{f_v}$, $m_{d_v}$ respectively represent the payload capacity, propellant (fuel) capacity, and dry mass.

The space mission planning problem ($P_0$ in Fig. \ref{fig:1b}) is different from the all-in-one formulation outlined in Section \ref{Problem Definition} with respect to the following two points: the nonlinear vehicle sizing constraint (Eq. \eqref{SL_vehicle_sizing}) is not included, and the quadratic penalty function is added to the objective function as Eq. \eqref{SubSL_obj} shows. Due to the quadratic objective function, this subproblem is a MIQP problem.

\begin{equation}
\begin{aligned}
\label{SubSL_obj}
\min_{\boldsymbol{x}_{vijt}, u_{vijt},\boldsymbol{y}_0} \quad &\sum_{t\in \mathcal{T}} \sum_{(v,i,j)\in \mathcal{A}} (\boldsymbol{a}_{vijt}^{T} \boldsymbol{x}_{vijt} + {a'}_{vijt}^{T} m_{d_{v}} u_{vijt}) + \phi_{0}\left(\boldsymbol{c}_{0}\left(\boldsymbol{y}, \boldsymbol{y}_{0}\right)\right) \\
\text{subject to} \quad & \text{Eqs.~\eqref{SL_constr1}--\eqref{SL_constr3} and \eqref{xdef}--\eqref{evdef} }\\
\text{where} \quad &
\boldsymbol{y} = [\boldsymbol{y}_1^T, \ldots , \boldsymbol{y}_N^T]^T \quad \text{and} \quad \boldsymbol{y}_v = [m_{p_v}, m_{f_v}, m_{d_v}]^T
\end{aligned}
\end{equation}

In the spacecraft design subproblems ($P_v$ in Fig. \ref{fig:1b}), the penalty function is minimized, and the vehicle sizing constraint ($m_{d_v} = \mathcal{F}(m_{p_v}, m_{f_v})$) is enforced. Note that this subproblem does not optimize the spacecraft design against its own cost or mass. As the objective of the original problem is to minimize IMLEO, the purpose here is to find a feasible spacecraft design minimizing the IMLEO value. The spacecraft design contains various interacting subsystems, and a hierarchical structure can be used to provide detailed subsystem-level design if needed. The subproblem for $v$-th type of vehicle can be expressed as Eq. \eqref{SizingSub}. Due to the nonlinear constraint, the subproblem is an NLP problem and can be solved by an NLP solver.

\begin{equation}
\begin{aligned}
\label{SizingSub}
\min_{\boldsymbol{y}_v} \quad &\phi_{v}\left(\boldsymbol{c}_{v}\left(\boldsymbol{y}, \boldsymbol{y}_{v}\right)\right) \\
\text{subject to} \quad &m_{d_v} = \mathcal{F}(m_{p_v}, m_{f_v})  \\
\text{where} \quad &
\boldsymbol{y}_v = [m_{p_v}, m_{f_v}, m_{d_v}]^T
\end{aligned}
\end{equation}

\subsection{Solution Algorithm and Iteration Scheme}
\label{iterscheme}
This subsection introduces the iterative solution algorithm for the decomposition-based algorithm introduced in Section \ref{OriginalALC} and Ref.~\cite{tosserams2007ALC}. The formulated decomposed optimization problems with ALC can be solved iteratively in two loops: the outer loop updates the augmented Lagrangian penalty parameters ($\boldsymbol{v}$, $\boldsymbol{w}$), while the inner loop solves the master problem and subproblems until the change in the subproblem objective function values is less than a defined tolerance. The outer loop terminates when two conditions are met simultaneously. First, the maximum consistency violation must be less than a specified tolerance. It ensures that the solution is feasible in each subproblem within the tolerance. Second, the maximum consistency violation change from the previous iteration must be smaller than the same tolerance. It allows for design space exploration and prevents premature convergences. The tolerance value of $10^{-3}$ is used in this study. The penalty parameters are updated in the outer loop based on the inner loop solution's consistency violation. Specifically, at $q$-th iteration, $\boldsymbol{v}\:  =\left[\boldsymbol{v_0}^{T}, \ldots, \boldsymbol{v}_{M-1}^{T}\right]^{T}$ is updated as follows:
\begin{equation}
    \boldsymbol{v}^{q+1}=\boldsymbol{v}^q+2\boldsymbol{w}^q\circ\boldsymbol{w}^q\circ\boldsymbol{c}^q
\end{equation}
The $r$-th element $w_r$ of the penalty weight vector $\boldsymbol{w}  =\left[\boldsymbol{w_0}^{T}, \ldots, \boldsymbol{w}_{M-1}^{T}\right]^{T}$ is updated based on the corresponding consistency constraint element $c_r$ as follows:
\begin{equation}
    w^{q+1}_r=
    \begin{cases}
    w^q_r &\text{if} \quad |c^q_r|\leq\gamma_2|c^{q-1}_r|\\
    \gamma_1 w^q_r &\text{if} \quad |c^q_r|>\gamma_2|c^{q-1}_r|
    \end{cases}
\end{equation}
\noindent where $\gamma_1>1$ and $0<\gamma_2<1$. In this study, $\gamma_1 = 2$ and $\gamma_2 = 0.5$ are used. The initial penalty parameter values can take $\boldsymbol{v}^1 = \bm{0}$ and $\boldsymbol{w}^1 \approx \bm{1}$.

The updates for the inner loop are performed by alternating between solving the master problem and the subproblems with the fixed penalty parameters. While each subproblem can be solved using the specialized numerical optimizer for MIQP or NLP, the master problem can be solved analytically as follows. 
\begin{equation}
\boldsymbol{y}= \underset{\boldsymbol{y}}{\operatorname{argmin}} \sum_{k=0}^{N} \phi_{k}\left(\boldsymbol{c}_{k}\left(\boldsymbol{y}, \boldsymbol{y}_{k}\right)\right)=\frac{\sum_{k=0}^{N}\left(\boldsymbol{w}_{k} \circ \boldsymbol{w}_{k} \circ \boldsymbol{y}_{k}\right)-\frac{1}{2} \sum_{k=0}^{N} \boldsymbol{v}_{k}}{\sum_{k=0}^{N}\left(\boldsymbol{w}_{k} \circ \boldsymbol{w}_{k}\right)}
\end{equation}

\begin{figure}[hbt!]
\centering
\includegraphics[width=.5\textwidth]{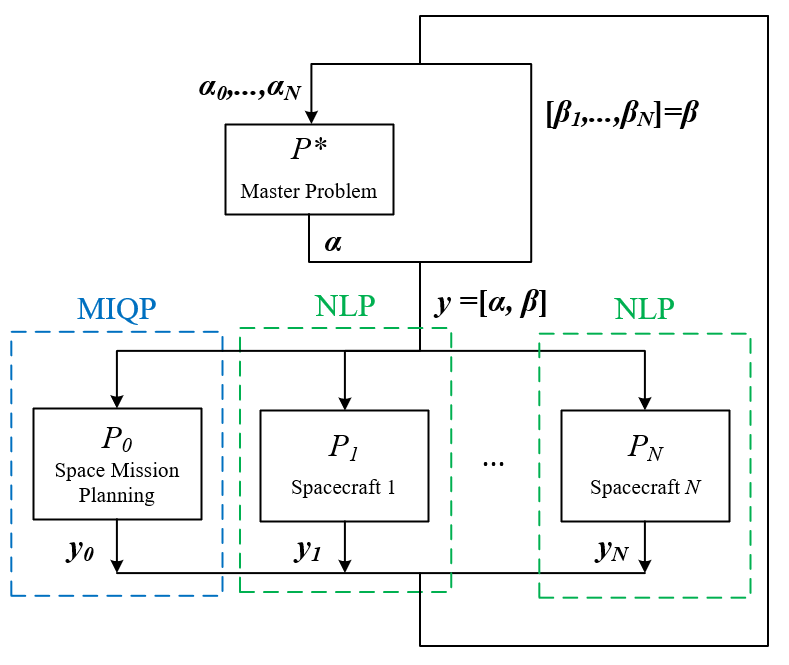}
\caption{Proposed decomposition-based optimization architecture with prioritized shared variables.}
\label{ALC_2}
\end{figure}

For our problem, we make an additional heuristics-based modification to the master problem to facilitate the convergence. The aforementioned master problem updates all the shared variables at the same time at every iteration, but this approach does not work effectively in our problem. This is because the space mission planning, with no knowledge of the constraints behind the spacecraft design, can return an aggressive or infeasible spacecraft design, which can deteriorate the convergence performance. Therefore, we propose to only update the spacecraft payload capacity and the propellant capacity in the master problem, while passing the spacecraft dry mass from the spacecraft design subproblem directly to the next iteration, as shown in Fig. \ref{ALC_2}. Mathematically, we separate the shared variables $\boldsymbol{y}$ into the regular shared variables $\boldsymbol{\alpha} = [m_{p_1}, m_{f_1}, \ldots, m_{p_N}, m_{f_N}]$ and the prioritized shared variables $\boldsymbol{\beta}=[m_{d_1}, \ldots, m_{d_N}]$ (i.e., $\boldsymbol{y} = [\boldsymbol{\alpha, \beta}]$), and only $\boldsymbol{\alpha}$ is updated in the master problem. 
\begin{equation}
\label{ModMaster}
\begin{aligned}
\min _ {\boldsymbol{\alpha}} \quad &\sum_{k=0}^{N} \phi_{k}\left(\boldsymbol{c}_{k}\left(\boldsymbol{\alpha}, \boldsymbol{\alpha}_k\right)\right)
\end{aligned}
\end{equation}
In the space mission planning subproblem, the spacecraft dry mass remains a variable, not a fixed parameter, and is subject to the penalty function. It indicates that the resultant dry mass $\beta_0$ is not used in the entire optimization architecture but only to facilitate the convergence of the whole optimization problem.

\subsection{Automatic Initial Solution Generation via Piecewise Linear Approximation}
\label{guess}
For the above iterative algorithm to perform effectively, a good initial guess of the shared variable is necessary. Thus, there is a need to develop an automatic and effective method that does not require a user-defined initial guess. To this end, we propose to approximate the nonlinear spacecraft model as a PWL model with chosen breakpoints so that an approximate initial solution to the original MINLP problem can be obtained using MILP \cite{chen2018MILP}. Since MILP problems can be solved using a specialized solver, this approximate solution can be generated efficiently. Although the PWL approximation does not necessarily return an optimal or even feasible solution to the original MINLP problem, the returned shared variables can be used as a good initial guess for the iterative approach.
Another advantage is that the MILP problem can be solved to the global optimum for the approximated nonlinear model. Thus, the MILP-based initial guess is not only automatically generated but also likely to be close to the nonlinear global optimum.

Specifically, in our problem, nonlinearity exists in the spacecraft sizing constraint. To generate the PWL functions of the nonlinear constraints, we choose a series of equally-spaced "mesh" points over the feasible ranges of the independent spacecraft design variables and use them as breakpoints for the PWL function generation. Since the dry mass is an (implicit) function of the payload capacity and propellant capacity, we only use the latter two for breakpoint generation. The breakpoint increment (or the number of breakpoints) is a key hyperparameter; a smaller increment (or more breakpoints) would lead to a more accurate initial guess, but it will also require a longer computational time. The method used in this study to model PWL functions as MILP problems is presented in Ref. \cite{PWLMILP}.

\section{Case Study: Human Lunar Exploration Campaign}
\label{CaseStudy}
To demonstrate the effectiveness of the proposed approach, we perform several case study instances and compare our approach with the state-of-the-art method. We first introduce the case study settings, followed by the results and the computational performance analysis.

\subsection{Case Study Settings}
\subsubsection{Specifications of Human Lunar Exploration Campaign}
\begin{figure}[hbt!]
\centering
\includegraphics[width=.6\textwidth]{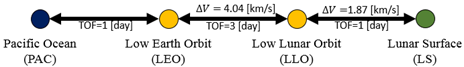}
\caption{Lunar campaign network model \cite{chen2018MILP}.}
\label{LunarMission}
\end{figure}

A human lunar exploration with two missions is considered here for the case study. 
The mission network model, parameters, commodity demand and supply used in this case study are presented in Fig. \ref{LunarMission}, Table \ref{MisisonParameter}, and Table \ref{LunarMissionDemand}, respectively. 
We call this case study the default instance, and seven additional instances with different parameters are also studied. Their specifications are given in Table \ref{CaseList}. The parameters not listed in Table \ref{CaseList} are the same as those in the default instance (Instance 1). In some instances, the supply and demand (i.e., mission payload) are increased, as shown in Table \ref{LunarMissionDemand}.
Note that sizing models for only one kind of spacecraft, which is a single-stage lunar lander, are considered for simplicity. Hence, the designed lander is also used as other kinds of spacecraft, such as in-space transfer vehicles. As landers are typically heavier than other spacecraft due to their landing structure, the optimization result might represent a conservative design. For this reason, a so-called \emph{aggressive} vehicle model with a lower structural mass estimation is used in some instances. In contrast, we refer to the original model as the conservative model. The sizing models are discussed in detail in the next section (\ref{sizing_ch}). For each case study instance, the computational time is measured on a platform with Intel Core i7-10700 (8 Core at 2.9 GHz). In the proposed decomposition-based method, Gurobi 9.1 solver \cite{gurobi} is used for the initial MILP problem and MIQP subproblem, and IPOPT \cite{IPOPT} is chosen for the NLP subproblems.

\begin{table}[h]
\caption{\label{MisisonParameter} Parameters used in the default instance (Instance 1) of the case study}
\centering
\begin{tabular}{p{0.5\textwidth}p{0.15\textwidth}}
\hline\hline Parameters & Assumed values \\\hline

Spacecraft Propellant type & LH2/LOX \\
Propellant $I_{sp}$, s & 420 \\
Propellant density $\rho_f$, kg/m\textsuperscript{3} & 360\\
Spacecraft miscellaneous mass fraction $c_{misc}$ (see Eq.~\eqref{sizing}) & 0.05\\
Spacecraft payload capacity range, kg & [500, 10000]\\
Spacecraft propellant capacity range, kg & [1000, 100000]\\
Type(s) of spacecraft designed  & 1 \\
Number of vehicles for each type & 6 \\
Crew mass (including space suit), kg/person & 100 \\
Crew consumption, kg/day/person & 8.655 \\
Spacecraft maintenance, structure mass/flight & 1\% \\

\hline\hline
\end{tabular}
\end{table}
\begin{table}[h]
\caption{\label{LunarMissionDemand} Lunar campaign commodity demand and supply}
\centering
\begin{tabular}{cccc}
\hline\hline Payload Type & Node & Date & Supply/Demand (increased instance) \\\hline
\multicolumn{4}{c}{Outbound to the Moon} \\\hline
Crew & Earth & 0, 365 & 4 \\
Habitat, Equipment, and Propellant, kg & Earth & 0, 365 & $\infty$ \\
Crew & Moon & 5, 370 & -4 \\
Habitat \& Equipment, kg & Moon & 5, 370 & -2000 (-3000)\\

\hline\multicolumn{4}{c}{Inbound to the Earth} \\\hline
Crew & Moon & 8, 373 & 4\\
Returned mass, kg & Moon & 8, 373 & 1000 (1500)\\
Crew & Earth & 13, 378 & -4 \\
Returned mass, kg & Earth & 13, 378 & -1000 (-1500)\\

\hline\hline
\end{tabular}
\end{table}
\begin{table}[hbt!]
\caption{\label{CaseList} Case study instances specifications}
\centering
\begin{tabular}{ccccc} 
\hline\hline
Instance index & \begin{tabular}[c]{@{}c@{}}Type(s) of SC\\ ~desinged\end{tabular} & \begin{tabular}[c]{@{}c@{}}Number of SC\\ per type\end{tabular} & SC sizing model & Mission payload  \\ 
\hline
1 (default)              & 1                                                                 & 6                                                               & conservative    & default          \\
2              & 2                                                                 & 3                                                               & conservative    & default          \\
3              & 6                                                                 & 1                                                               & conservative    & default          \\
4              & 2                                                                 & 3                                                               & conservative    & increased        \\
5              & 1                                                                 & 6                                                               & aggressive      & default          \\
6              & 2                                                                 & 3                                                               & aggressive      & default          \\
7              & 6                                                                 & 1                                                               & aggressive      & default          \\
8              & 2                                                                 & 3                                                               & aggressive      & increased        \\
\hline\hline
\end{tabular}
\end{table}

\clearpage
\subsubsection{Spacecraft Sizing Model}\label{sizing_ch}
The subsystem-level spacecraft model used as the spacecraft design constraint in Eq. \eqref{SL_vehicle_sizing} is developed by the least square curve fitting to the data from the lunar lander design database in Ref.~\cite{isaji2018landerdata, isaji2020lander}. The following set of equations shows the spacecraft model used in the case study.
\begin{align}
\begin{aligned}
\label{sizing}
 &m_d=\sum m_{sub} = m_{str}+m_{prop}+m_{power}+m_{avi}+m_{ECLSS}+m_{misc}
    \\[0.5em] & \text { where }\\
    & \begin{array}{ll}
        \quad m_{str} &= {n_{stg}}^{-0.6705}\,(0.3238\,{m_{d}} + 693.7\,{m_{p}}^{0.04590}) 
        \\[1em]
        \quad m_{prop} &= {0.1648}\,(m_{d}+m_{p})+20.26\,\left(\dfrac{m_{f}}{\rho_{f}}\right)
        \\[1em]
        \quad m_{power} &= 7.277\cdot10^{-8}\,{m_{d}}^{2.443}+137.0
        \\[1em]
        \quad m_{avi} &= 1.014\,{m_{power}}^{0.8423}+22.33\,{t_{mis}}
        \\[1em]
        \quad m_{ECLSS} &= 0.004190\,n_{crew}\,t_{mis}\,{m_{d}}^{0.9061}\,{n_{stg}}^{0.7359}+434.7
        \\[1em]
        \quad m_{misc} &= c_{misc}\, m_{d}
    \end{array}
\end{aligned}
\end{align}

Beyond the payload capacity and propellant capacity, there are some additional parameters in these equations: $n_{stg}$ is the number of stages (either 1 or 2), $\rho_f$ is the propellant density in kg/m\textsuperscript{3}, $t_{mis}$ is the surface time of the lunar mission in days, $n_{crew}$ is the number of crew, $c_{misc}$ is the miscellaneous mass fraction. The miscellaneous mass fraction $c_{misc}$ represents how much of the dry mass is categorized as the miscellaneous mass. It can range from 0 to 0.15, meaning 0\% to 15\% of the dry mass is the miscellaneous mass. The higher $c_{misc}$, the heavier and more conservative the vehicle design becomes. All mass properties are defined in kg.

As shown in Eq.~\eqref{sizing}, the model captures the subsystem-level interactions to return the relationship between the payload capacity, propellant capacity, and dry mass of the spacecraft. Particularly, the subsystem interactions are captured \emph{through} the dry mass. For instance, an increase in any subsystem mass will raise the dry mass. Since each subsystem mass is dependent on the dry mass, their mass should increase as well, which will further raise the dry mass. The 'balanced' dry mass with such subsystem circular references can be found by solving Eq. \eqref{sizing} for the dry mass, $m_d$. More details on this model can be found in Appendix A.

The aforementioned aggressive sizing model has a lower structural mass estimation shown in Eq. \eqref{aggressive}. The only difference from Eq. \eqref{sizing} is the lowered coefficient of $m_d$, and the other subsystem mass relations remain the same. In the ranges of payload capacity and propellant capacity specified in Table \ref{MisisonParameter}, both models are non-convex.

\begin{equation}
\label{aggressive}
    m_{str} = {n_{stg}}^{-0.6705}\,(0.2694\,{m_{d}} + 693.7\,{m_{p}}^{0.04590})
\end{equation}

\subsection{Optimization Results of the Proposed Decomposition-Based Formulation}
\label{ALC_CS}
This subsection introduces the optimization results of the proposed decomposition-based formulation. Since the performance of the proposed method is dependent on the breakpoint increment (or the number of breakpoints) for the PWL approximation of the MILP-based initial solution generation, five different increments are tested for Instance 1 \footnote{We have tested the integrated formulation (as a black box) with MINLP solvers including Outer Approximation \cite{OuterApproximation}, LP/NLP based Branch-and-Bound \cite{LPNLP-BB}, and Extended Cutting Plane \cite{ExtendedCuttingPlane}; however, they failed to produce a feasible solution due to the non-convex sizing constraints.}. The results are shown in Table \ref{ALCdefault}; note that the PWL-based initial solutions are not necessarily feasible, and so their IMLEOs are only reported for information purposes.
\begin{table}[!htbp]
\centering
\caption{\label{ALCdefault}Instance 1 optimization results by the proposed decomposition-based formulation}
\begin{tabular}{ccccclcc} 
\hline\hline
\multirow{3}{*}{\begin{tabular}[c]{@{}c@{}}PWL\\increment, kg\end{tabular}} & \multirow{3}{*}{\begin{tabular}[c]{@{}c@{}}PWL~ \\mesh points\end{tabular}} & \multicolumn{3}{c}{Optimization time, s}                             &                      & \multicolumn{2}{c}{IMLEO, kg}                \\ 
\cline{3-5}\cline{7-8}
                                                                            &                                                                             & PLW-based            & Decomposition-based        & \multirow{2}{*}{Total} & \multicolumn{1}{c}{} & PWL-based            & Final                 \\
                                                                            &                                                                             & initial solution     & iterations         &                        & \multicolumn{1}{c}{} & initial solution     & solution              \\ 
\hline
10,000                                                                      & 13                                                                          & 3.461                & 20.02                & 23.48                  &                      & 741,115              & 724,776               \\
5,000                                                                       & 36                                                                          & 4.769                & 13.47                & 18.24                  &                      & 700,684              & 694,224               \\
2,500                                                                       & 120                                                                         & 4.754                & 11.98                & 16.73                  &                      & 677,035              & 676,862               \\
1,250                                                                       & 425                                                                         & 28.09                & 13.83                & 41.92                  &                      & 677,343              & 677,204               \\
625                                                                         & 1,595                                                                       & 465.5                & 13.18                & 478.7                  &                      & 677,315              & 677,072               \\ 
\hline\hline
\multicolumn{1}{l}{}                                                        & \multicolumn{1}{l}{}                                                        & \multicolumn{1}{l}{} & \multicolumn{1}{l}{} & \multicolumn{1}{l}{}   &                      & \multicolumn{1}{l}{} & \multicolumn{1}{l}{} 
\end{tabular}
\end{table}

Although the optimizer's computational time involves some randomness depending on the individual problems, some general trends can be observed. First, when the increment is too large (too few breakpoints, e.g., 10,000 kg increment with 13 mesh points), the initial solution quality is poor, and thus the final solution IMLEO is also poor. 
Second, the computational time to solve the initial MILP problem rapidly increases when the increment is too small (too many breakpoints, e.g., 625 kg increment with 1,595 mesh points), resulting in a long total computational time.
In summary, we can observe the expected trend that a smaller increment (more breakpoints) leads to a better initial guess at the cost of computational time. 
Theoretically speaking, if we reduce the increment to zero (an infinite number of breakpoints), the solution would match with the global optimum; however, this is impractical as it requires infinite computational time. 
Thus, the most efficient strategy is to use an increment that can generate a reasonably accurate initial solution and leave the rest to the proposed decomposition-based optimization. 
Although this hyperparameter needs to be chosen for the proposed algorithm, it is worth noting that the computational performance is not very sensitive to the choice of its value except for the extreme cases. In our case study and spacecraft sizing model, the 2,500 kg increment is chosen as the default increment for the remainder of the case study.

With the chosen default increment, the other seven instances are also tested. The results are summarized in Table \ref{ALCinstances}. In all instances, the proposed method can solve the problem with a reasonable amount of computational time and solution improvement from the initial guess. Note that the final IMLEO is not necessarily lower than that of the PWL-based initial solution because the original nonlinear problem might be infeasible with the initial guess spacecraft design (i.e., $y_v^0 = [m_{p_v}^0, m_{f_v}^0, \mathcal{F}(m_{p_v}^0, m_{f_v}^0)]^T$ for $v$-th vehicle).
\begin{table}[!htbp]
\centering
\caption{\label{ALCinstances}Decomposition-based optimization results with 2,500 kg increment for all instances}
\begin{tabular}{cccclcc} 
\hline\hline
            & \multicolumn{3}{c}{Optimization time, s}                        &                      & \multicolumn{2}{c}{IMLEO, kg}  \\ 
\cline{2-4}\cline{6-7}
Instance    & PWL-based        & Decomposition-based & \multirow{2}{*}{Total} & \multicolumn{1}{c}{} & PWL-based        & Final       \\
            & initial solution & iterations        &                        & \multicolumn{1}{c}{} & initial solution & solution    \\ 
\hline
1 (default) & 4.754            & 11.98               & 16.73                  &                      & 677,035          & 676,862     \\
2           & 38.66            & 10.70               & 49.35                  &                      & 401,191          & 401,093     \\
3           & 38.66            & 61.57               & 100.2                  &                      & 386,257          & 387,535     \\
4           & 10.53            & 22.63               & 33.16                  &                      & 470,788          & 470,406     \\
5           & 75.35            & 94.06               & 154.7                  &                      & 442,596          & 442,605     \\
6           & 24.09            & 14.30               & 38.39                  &                      & 293,095          & 293,095     \\
7           & 100.1            & 44.02               & 144.1                  &                      & 286,154          & 302,041     \\
8           & 31.50            & 14.47               & 45.97                  &                      & 344,397          & 344,423     \\
\hline\hline
\end{tabular}
\end{table}

Overall, the proposed decomposition-based formulation can take the reasonable approximate solution by the PWL formulation and offer a better computational efficiency to achieve a high-quality solution. 

\subsection{Benchmark Formulation: Modified Embedded Optimization}
\label{embedded_CS}
While our formulation of the integrated mission planning and subsystem-level spacecraft design has not been directly tackled in the literature, we can extend a state-of-the-art approach for a similar problem straightforwardly as a benchmark to evaluate our newly proposed method. The identified approach is the embedded optimization method by Taylor \cite{taylor2007phd}, which was demonstrated to be more efficient than directly solving the original integrated MINLP problem using a global optimizer. In this paper, the method is modified to solve the case study problem more efficiently.

In the original formulation of the embedded optimization, all spacecraft variables are separated from the problem and determined by a metaheuristics algorithm. However, in our problem, the spacecraft dry mass is constrained and can be uniquely determined by the payload and propellant capacity. Hence, we let the metaheuristics algorithm pick only the payload and propellant capacity of $N$ types of vehicles, and calculate the corresponding spacecraft dry mass by the spacecraft sizing constraint. This process reduces the number of variables from 3$N$ to 2$N$ and turns the constrained problem into an unconstrained one.
After obtaining the feasible vehicle design, these values are fed to the space mission planning problem, which is solved by the specialized MILP optimizer. Unlike the all-in-one formulation, the vehicle parameters are fixed within the mission planning part. Then, the corresponding objective function value (i.e., IMLEO) is returned to the metaheuristic optimizer for the evaluation for the next iteration.
As a result, the metaheuristics only handles an unconstrained optimization problem with $2N$ variables (i.e., the payload capacity and propellant capacity for $N$ spacecraft), where the evaluation of the constraints and the determination of the remaining variables are handled by the embedded MILP solver. When the problem is infeasible with the chosen spacecraft design, a predetermined large positive number is returned as the objective function value, following the death-penalty constraint handling \cite{DeathPenalty}; this method is chosen because the large number of constraints in the problem makes other methods (e.g., repair-based methods) challenging or inapplicable. The problem to be optimized by the metaheuristic solver is expressed as Eq. \eqref{heuristics}.

\begin{equation}
\label{heuristics}
\begin{aligned}
\min_{\boldsymbol{\alpha}} \quad &\text{IMLEO}(\boldsymbol{\alpha}, \mathcal{F}(\boldsymbol{\alpha}))\\
\text{where} \quad &\boldsymbol{\alpha} = [m_{p_1}, m_{f_1}, \ldots, m_{p_N}, m_{f_N}], \quad\boldsymbol{\alpha} \in \mathbb{R}^{2N}
\end{aligned}
\end{equation}

In accordance with the proposed decomposition-based method, the PWL MILP initial guess is given to the initial population.  Since such an initial guess is not provided in the original method, this modification should lead to performance improvement. For the same reason as the proposed method, a reasonable increment for the PWL approximation is desired. The effect of the PWL increment value is examined (see Section \ref{increment&algo}). 

Since the performance of the embedded optimization would depend on the choice of the metaheuristics algorithm, three different metaheuristics algorithms are tested: the extended Ant Colony Optimization (ACO) \cite{AntColony}, the Genetic Algorithm (GA) \cite{GA}, and the Particle Swarm Optimization (PSO) \cite{PSO}. The population size is 10 for all, and other specific parameters to each algorithm are provided in Appendix B. The optimization is terminated when a predefined number of generations are populated; progress at different generation numbers are examined for each algorithm to explore the tradeoff between the computational time and accuracy. Due to the random nature of the metaheuristic optimizers, the optimization is run 10 times with the same setting. 

\subsection{Performance Comparison of the Proposed Formulation with Modified Embedded Optimization}
In this section, the proposed formulation results are compared with the modified embedded optimization results, which depend on several factors, including the PWL increment (i.e., initial guess quality), algorithm, and problem to be solved. We first compare their performance with various increments and algorithms, then study the results of all case study instances.
A set of sample results are shown in this section to discuss the findings; the complete dataset for every instance and algorithm is presented in Appendix B.
\subsubsection{Performance Comparison with Various PWL Increments}\label{increment&algo}

First, we examine the results for the default instance with different increments. Table \ref{10gen} shows how the three metaheuristics algorithms compare with the proposed method at 10 generations for various PWL increments. A generation number as low as 10 is used to match the time scale of the proposed method. Note that the optimization time in Table \ref{10gen} does not include the initial guess generation time because it is the same as long as the increment is the same.
\begin{table}[!t]
\centering
\caption{\label{10gen}Instance 1 results comparison with different PWL increments at 10 generations}
\begin{tabular}{ccccc} 
\hline\hline
PWL increment, kg       & Algorithm & Average opt. time, s & Best IMLEO, kg & Worst IMLEO, kg  \\ 
\hline
\multirow{4}{*}{10,000} & ACO       & 66.93        & 707,882        & 724,776          \\
                        & GA        & 72.81        & 684,755        & 724,776          \\
                        & PSO       & 69.06        & 696,625        & 724,776          \\
                        & Proposed  & 20.02        & 724,776        & -                \\ 
\hline
\multirow{4}{*}{5,000}  & ACO       & 65.99        & 694,224        & 694,224          \\
                        & GA        & 71.27        & 681,842        & 694,224          \\
                        & PSO       & 67.98        & 694,224        & 694,224          \\
                        & Proposed  & 13.47        & 694,224        & -                \\ 
\hline
\multirow{4}{*}{2,500}  & ACO       & 72.27        & 705,550        & infeasible       \\
                        & GA        & 80.58        & 759,326        & infeasible       \\
                        & PSO       & 74.57        & 700,227        & infeasible       \\
                        & Proposed  & 11.98        & 676,862        & -                \\ 
\hline
\multirow{4}{*}{1,250}  & ACO       & 65.93        & 677,303        & 677,303          \\
                        & GA        & 70.16        & 677,303        & 677,303          \\
                        & PSO       & 69.00        & 677,303        & 677,303          \\
                        & Proposed  & 13.83        & 677,204        & -                \\ 
\hline
\multirow{4}{*}{625}    & ACO       & 67.21        & 677,283        & 677,283          \\
                        & GA        & 69.42        & 677,283        & 677,283          \\
                        & PSO       & 68.33        & 677,283        & 677,283          \\
                        & Proposed  & 13.18        & 677,072        & -                \\
\hline\hline
\end{tabular}
\end{table}

There are several key findings. Most importantly, the proposed method converges substantially faster than any metaheuristics algorithms, even if the number of generations for metaheuristics is limited to 10. In most cases, the converged solution of the proposed method is better than even the best case (out of the 10 runs) for the modified embedded optimization. Note that it is not surprising that metaheuristics can occasionally perform better than the proposed method because it can stochastically escape from the local minima; however, it is worth stressing that, due to the stochasticity, the performance of the embedded optimization can vary significantly and possibly return infeasible solution in the worst case; this is true for the increment of 2,500 kg, where the initial PWL-based solution is not feasible and the embedded optimization is not able to lead it to a feasible one in the worst case of the 10 runs. In contrast, the performance of the proposed method is robust and repeatable due to its deterministic nature; it can find a nonlinear (local) solution regardless of the initial guess feasibility even though it might not always improve significantly from the initial guess.

\subsubsection{Performance Comparison in Different Problem Instances}\label{inst}
Next, we select PSO as our default algorithm and compare the results for all eight instances. 
The results by PSO for all instances are shown in Table \ref{PSOinstances}. A similar trend is observed for all cases: the results at 10 generations suggest that there is a robust trend that the proposed method can achieve a better solution given a similar amount of time. 
In fact, note that in all cases, the computational time of the proposed method is substantially shorter than any other methods. 
\begin{table}[!b]
\centering
\caption{\label{PSOinstances}PSO embedded opt. results comparison: 2,500 kg PWL increment, all instances, 10 and 100 generations}
\begin{tabular}{cccccc} 
\hline\hline
Instance                     & Algorithm            & Number of generations & Average opt. time, s & Best IMLEO, kg & Worst IMLEO, kg  \\ 
\hline
\multirow{3}{*}{1 (default)} & \multirow{2}{*}{PSO} & 10                    & 74.57        & 700,227        & infeasible       \\
                             &                      & 100                   & 765.3        & 680,234        & 710,650          \\
                             & Proposed             & -                     & 11.98        & 676,862        & -                \\ 
\hline
\multirow{3}{*}{2}           & \multirow{2}{*}{PSO} & 10                    & 84.70        & 436,363        & infeasible       \\
                             &                      & 100                   & 881.6        & 402,976        & 412,622          \\
                             & Proposed             & -                     & 10.70        & 401,093        & -                \\ 
\hline
\multirow{3}{*}{3}           & \multirow{2}{*}{PSO} & 10                    & 153.0        & 423,873        & infeasible       \\
                             &                      & 100                   & 1620         & 408,388        & 421,829          \\
                             & Proposed             & -                     & 61.57        & 387,535        & -                \\ 
\hline
\multirow{3}{*}{4}           & \multirow{2}{*}{PSO} & 10                    & 82.52        & 470,538        & 470,538          \\
                             &                      & 100                   & 842.3        & 470,538        & 470,538          \\
                             & Proposed             & -                     & 22.63        & 470,406        & -                \\ 
\hline
\multirow{3}{*}{5}           & \multirow{2}{*}{PSO} & 10                    & 75.36        & 451,269        & 470,869          \\
                             &                      & 100                   & 748.6        & 445,230        & 450,841          \\
                             & Proposed             & -                     & 94.06        & 442,605        & -                \\ 
\hline
\multirow{3}{*}{6}           & \multirow{2}{*}{PSO} & 10                    & 96.64        & 305,457        & 337,105          \\
                             &                      & 100                   & 949.7        & 298,414        & 306,893          \\
                             & Proposed             & -                     & 14.30        & 293,095        & -                \\ 
\hline
\multirow{3}{*}{7}           & \multirow{2}{*}{PSO} & 10                    & 172.3        & 302,791        & 358,768          \\
                             &                      & 100                   & 1772         & 298,361        & 310,755          \\
                             & Proposed             & -                     & 44.02        & 302,041        & -                \\ 
\hline
\multirow{3}{*}{8}           & \multirow{2}{*}{PSO} & 10                    & 94.21        & 355,767        & 395,748          \\
                             &                      & 100                   & 934.3        & 348,585        & 359,040          \\
                             & Proposed             & -                     & 14.47        & 344,423        & -                \\
\hline\hline
\end{tabular}
\end{table}

It is of interest to observe the performance of the metaheuristic methods with a higher number of generations. Thus, the 100-generation case is also shown in Table \ref{PSOinstances}. In this case, we can see that there are (rare) cases where PSO can find a better solution than the proposed method in the best case (e.g., Instance 7); however, as mentioned previously, the stochastic nature of the metaheuristics algorithms makes their performance highly variable; for all cases tested, the proposed method reached the same or better solutions than those from the worst-case metaheuristics method even after 100 generations and with a substantially shorter computational time.

\subsection{Advantages and Limitations of the Proposed Formulation}
One advantage of the proposed method is its potential to integrate more complicated spacecraft design problems, such as models with more constraints or even MDO problems. Oftentimes, MDO problems require external and black-box analysis tools for specific disciplines where the analytical formulation or constraints might not be known by the user. As the ALC formulation is intended for such MDO problems, integration of disciplinary analysis tools should not be an issue. The method we employ to model PWL functions as MILP problems can handle black-box functions as well \cite{PWLMILP}. Furthermore, if a hierarchical MDO subproblem is considered, a multi-level hierarchical variation of the ALC formulation can be utilized \cite{ATCextended}. When there are a large number of subproblems involved, parallel computation might further reduce the computational cost.

The proposed method has several limitations as well, which can lead to future work. First, the objective function of the space mission planning problem must be linear or quadratic. Otherwise, the subproblem itself becomes a MINLP problem and negates the advantage of the proposed decomposition formulation. This limitation excludes general nonlinear objective functions. In a similar manner, the spacecraft design constraints can only have continuous variables although integer variables such as the number of engines and number of propellant tanks may be needed for high-fidelity spacecraft design. One approach to these cases is to revise the decomposition structure of the problem, which is left for future work. Finally, even though the proposed PWL solution can provide a near-optimal initial guess, the gradient-based NLP solver in the proposed method can still lead to a local optimum. Thus, embedded optimization may be preferred in certain contexts when the goal is to achieve a global optimum without any constraint on the available computational resources.

Overall, the case study demonstrates the high computational performance of the proposed method. The proposed formulation can achieve a high-quality solution robustly in a shorter computational time than the state-of-the-art modified embedded optimization method. It is also illustrated that the computational efficiency is not impaired by infeasible initial solutions.

\section{Conclusion}
\label{conclusion}
This paper tackles the challenging problem of integrated space mission planning and spacecraft design. The all-in-one formulation is presented as an MINLP problem, and an efficient solution approach is developed leveraging the unique structure of the problem and following the philosophy of MDO. Namely, the all-in-one MINLP problem is decomposed into the space mission planning subproblem (MIQP) and the spacecraft design subproblem(s) (NLP) so that they can be solved iteratively using the ALC approach to find the optimal solution for the original MINLP problem. Furthermore, an automatic and effective approach for finding an initial solution for this iterative process is proposed leveraging a piecewise linear (PWL) approximation of the nonlinear vehicle model, so that no user-defined initial guess is needed.
In the case study, the computational efficiency of the proposed method is studied with several problem settings.
The results demonstrate that, compared to the state-of-the-art method, the proposed formulation converges significantly faster, and the converged
solution is also at least the same or better in the same computational time limit, even when started from the same PWL-based initial solution. The combination of the unique problem structure, the iterative algorithms for shared variables, and the efficient initial solution generation method leads to this computational efficiency. The parallelizable nature of the algorithm can potentially make the proposed method even more effective for large-scale problems. Due to the flexibility of the ALC method, the proposed formulation can also integrate more complex vehicle design models, which is left for future work.

\section*{Appendix A: Spacecraft Design Model}
\label{App_A}
\begin{figure}[hbt!]
\centering
\includegraphics[width=.45\textwidth]{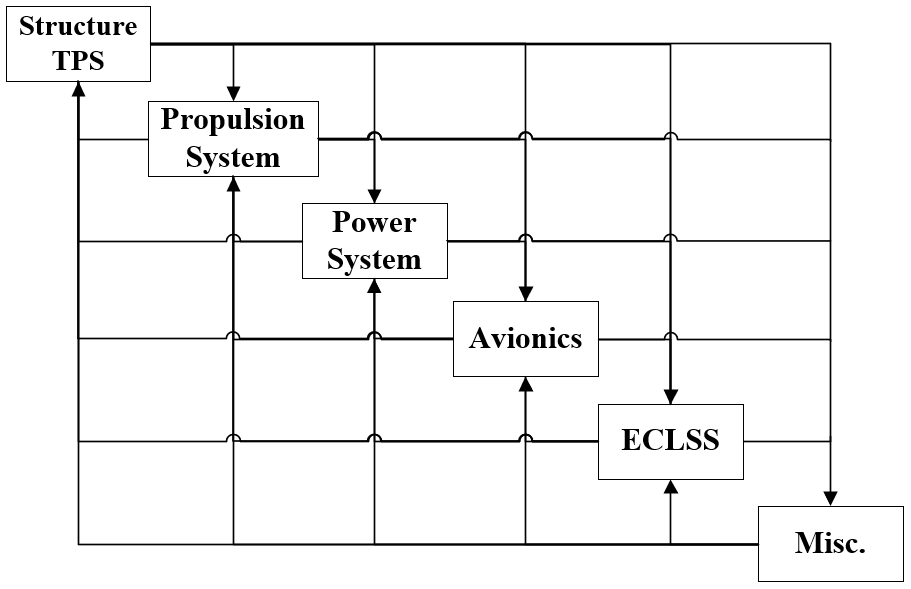}
\caption{Relationship of domains in a single-stage lunar lander.}
\label{N2}
\end{figure}

This appendix provides more details on the parametric sizing model for the spacecraft used in the case study. 
In the considered model, the subsystems of single-stage landers and their relations to the dry mass are defined as Eq. \eqref{drymass_eq}. 

\begin{equation}
\label{drymass_eq}
m_{d}=\sum m_{sub} = m_{str}+m_{prop}+m_{power}+m_{avi}+m_{ECLSS}+m_{misc}
\end{equation}

\noindent where $m_{sub}$ indicates the mass of each subsystem. $m_{str}$ indicates the structure and thermal protection system (TPS), which includes all subsystems that are attached to support or connect other components. This includes landing legs and truss, TPS for the reentry to the earth, and a docking mechanism. $m_{prop}$ is the propulsion system, such as propellant tanks, reaction control system (RCS), and hardware of engines. $m_{power}$ is the power system, which contains batteries, fuel cells, solar panels, or other electrical systems. $m_{avi}$ indicates the avionics, and $m_{ECLSS}$ indicates the environmental and life control system (ECLSS) that supports the crew's lives such as consumables (food, water, air) or related piping and tankage. Finally, we also consider other miscellaneous required components, expressed as $m_{misc}$. Through the dry mass, each subsystem interacts with every other subsystem, and this relation is visualized in Fig.~\ref{N2} as an N\textsuperscript{2} diagram.

For the defined subsystems, mass estimation relationships (MERs) are developed as functions of payload capacity, propellant capacity, propellant type $\zeta$, and some other known parameters. If the propellant type is fixed, the subsystems MERs and dry mass are dependent on the payload capacity and propellant capacity only, and thus serve as the vehicle sizing constraints (Eq. \eqref{SL_vehicle_sizing}, $m_d = \mathcal{F}(m_p, m_f$)). Each subsystem MER is developed by the least square curve fitting to the data from the lunar lander design database in Ref.~\cite{isaji2018landerdata, isaji2020lander}, which includes both existing and elaborated conceptual design.
The form of each subsystem's MER is manually determined to be a sufficiently simple yet accurate form. 
The resultant MERs are shown in Eq. \eqref{sizing}.

Table \ref{tab_SizingVar} summarizes the independent variables, $R^2$ values for curve fitting, number of data points used for curve fitting ($N_{data}$), average errors against the data points, and the maximum errors. Note that only a small number of data points are used for the propulsion system MER since two-stage lander data are excluded as their propulsion systems with staging are too distinct from those of single-stage ones. One may also see that relatively poor correlations are obtained for the power systems and avionics mass as they simply might not be strong functions of the dry mass or vehicle size. However, since they typically account for small portions of the dry mass, the poor correlation does not have a significant effect on the validation process. 

The limitation of this sizing model should also be noted. Because the MERs are developed from the existing data points, a solution for vehicles that are significantly heavier than the ones in the database would either be a low-fidelity model or infeasible. In other words, $m_d$ that satisfies Eq. \eqref{drymass_eq} might not exist for certain weight classes. Specifically, the upper bound of the dry mass is approximately 23,000 kg. When $t_{mis}$ is 3 days, $n_{crew}$ is 4, $c_{misc}$ is 0.05, and the propellant is LH2/LOX, the upper bound are found at 500 kg payload and 75,500 kg propellant, or at 10,000 kg payload and 45,500 kg propellant.

\begin{table}[h]
\caption{\label{tab_SizingVar} Summary of subsystem MERs }
\centering
\begin{tabular}{p{0.17\textwidth}p{0.09\textwidth}p{0.2\textwidth}p{0.08\textwidth}p{0.08\textwidth}p{0.1\textwidth}p{0.1\textwidth}}
\hline\hline Subsystem & Notation & Independent Variables & $R^2$ & $N_{data}$ & Avg. Error & Max. Error\\
\hline
Structure + TPS & $m_{str}$ & $m_{d},n_{stg},m_p$ & 0.9254 & 17 & 7.379\% & 24.31\% \\
Propulsion System & $m_{prop}$ & $m_d, m_p, \rho_p$ & 0.9279 & 8 & 7.429\% & 11.16\%\\
Power System & $m_{power}$ & $m_d$ & 0.7182 & 13 & 16.24\% & 36.68\%\\
Avionics & $m_{avi}$ & $m_{power}(m_d), t_{mis}$ & 0.6204 & 22 & 36.42\% & 75.94\%\\
ECLSS & $m_{ECLSS}$ & $m_d, n_{crew}, n_{stg}, t_{mis}$ & 0.9293 & 12 & 11.93\% & 38.09\% \\
Miscellaneous & $m_{misc}$ & $m_d$ & - & - & - & - \\
\hline\hline
\end{tabular}
\end{table}

\section*{Appendix B: Summary of the Embedded Optimization Results}
Table \ref{ACOall}, \ref{GAall}, and \ref{PSOall} show the modified embedded optimization results by ACO, GA, and PSO, respectively. The initial guess increment is 2,500 kg for all cases. The algorithm implementations by pygmo \cite{pygmo}, a scientific Python library, are used to generate the data. The parameters used for each algorithm are as follows.

For ACO, the kernel size, convergence speed parameter, oracle parameter, accuracy parameter, threshold parameter, and focus parameter are 10, 1, $10^9$, 0, 7, 0, respectively.
GA uses the exponential crossover \cite{DE} with 0.9 probability, polynomial mutation \cite{PolMutation} with 0.02 probability, and tournament selection with size 2.
For PSO, the inertia weight, social component, cognitive component, and maximum particle velocities are 0.7298, 1.05, 2.05, and 0.5, respectively.

\begin{table}[hbt!]
\centering
\caption{\label{ACOall}ACO embedded optimization results with 2,500 kg PWL increment}
\begin{tabular}{ccccc} 
\hline\hline
Instance & Number of generations & Optimization time, s & \multicolumn{1}{l}{Best IMLEO, kg} & \multicolumn{1}{l}{Worst IMLEO, kg}  \\ 
\hline
         & 10                    & 72.27                & 705,550                            & infeasible                                  \\
1        & 50                    & 360.9                & 705,550                            & 830,350                              \\
         & 100                   & 721.2                & 685,146                            & 783,343                              \\ 
\hline
         & 10                    & 91.37                & 461,323                            & 831,125                              \\
2        & 50                    & 452.4                & 420,662                            & 631,469                              \\
         & 100                   & 904.0                & 420,662                            & 541,799                              \\ 
\hline
         & 10                    & 147.0                & 460,524                            & infeasible                                  \\
3        & 50                    & 732.0                & 459,157                            & 565,537                              \\
         & 100                   & 1466                 & 426,322                            & 508,066                              \\ 
\hline
         & 10                    & 80.68                & 470,538                            & 470,538                              \\
4        & 50                    & 403.6                & 470,538                            & 470,538                              \\
         & 100                   & 804.6                & 470,538                            & 470,538                              \\ 
\hline
         & 10                    & 74.45                & 454,849                            & 484,095                              \\
5        & 50                    & 371.3                & 453,073                            & 466,563                              \\
         & 100                   & 740.9                & 445,776                            & 461,279                              \\ 
\hline
         & 10                    & 95.43                & 321,152                            & 359,303                              \\
6        & 50                    & 476.7                & 299,690                            & 336,976                              \\
         & 100                   & 952.5                & 299,690                            & 330,914                              \\ 
\hline
         & 10                    & 171.8                & 323,436                            & 348,874                              \\
7        & 50                    & 858.8                & 315,515                            & 326,136                              \\
         & 100                   & 1716                 & 313,764                            & 323,436                              \\ 
\hline
         & 10                    & 94.14                & 375,746                            & 425,749                              \\
8        & 50                    & 468.6                & 359,239                            & 401,878                              \\
         & 100                   & 934.7                & 351,292                            & 394,226                              \\
\hline\hline
\end{tabular}
\end{table}
\begin{table}[hbt!]
\centering
\caption{\label{GAall}GA embedded optimization results with 2,500 kg PWL increment}
\begin{tabular}{ccccc} 
\hline\hline
Instance & Number of generations & Optimization time, s & \multicolumn{1}{l}{Best IMLEO, kg} & \multicolumn{1}{l}{Worst IMLEO, kg}  \\ 
\hline
         & 10                    & 80.58                & 759,326                            & infeasible                                  \\
1        & 50                    & 413.1                & 712,515                            & infeasible                                  \\
         & 100                   & 814.8                & 712,515                            & infeasible                                  \\ 
\hline
         & 10                    & 98.80                & 437,974                            & infeasible                                  \\
2        & 50                    & 511.9                & 437,275                            & infeasible                                  \\
         & 100                   & 1010                 & 414,034                            & infeasible                                  \\ 
\hline
         & 10                    & 158.7                & 427,682                            & infeasible                                  \\
3        & 50                    & 840.4                & 419,610                            & infeasible                                  \\
         & 100                   & 1686                 & 398,507                            & infeasible                                  \\ 
\hline
         & 10                    & 85.71                & 470,538                            & 470,538                              \\
4        & 50                    & 434.4                & 470,538                            & 470,538                              \\
         & 100                   & 872.8                & 470,538                            & 470,538                              \\ 
\hline
         & 10                    & 84.28                & 455,676                            & 484,630                              \\
5        & 50                    & 406.7                & 455,676                            & 479,622                              \\
         & 100                   & 806.1                & 455,676                            & 476,852                              \\ 
\hline
         & 10                    & 106.2                & 306,077                            & 359,839                              \\
6        & 50                    & 512.0                & 297,363                            & 329,513                              \\
         & 100                   & 1020                 & 297,183                            & 316,152                              \\ 
\hline
         & 10                    & 178.5                & 304,792                            & 345,599                              \\
7        & 50                    & 866.8                & 296,631                            & 323,739                              \\
         & 100                   & 1708                 & 291,836                            & 305,488                              \\ 
\hline
         & 10                    & 111.7                & 358,082                            & 513,027                              \\
8        & 50                    & 518.8                & 353,113                            & 388,664                              \\
         & 100                   & 998.2                & 349,502                            & 371,646                              \\
\hline\hline
\end{tabular}
\end{table}
\begin{table}[!htbp]
\centering
\caption{\label{PSOall}PSO embedded optimization results with 2,500 kg PWL increment}
\begin{tabular}{ccccc}
\hline\hline
Instance & Number of generations & Optimization time, s & \multicolumn{1}{l}{Best IMLEO, kg} & \multicolumn{1}{l}{Worst IMLEO, kg}  \\ 
\hline
         & 10                    & 74.57                & 700,227                            & infeasible                                  \\
1        & 50                    & 381.7                & 680,234                            & 731,419                              \\
         & 100                   & 765.3                & 680,234                            & 710,650                              \\ 
\hline
         & 10                    & 84.70                & 436,363                            & infeasible                                  \\
2        & 50                    & 440.0                & 403,448                            & 439,154                              \\
         & 100                   & 881.6                & 402,976                            & 412,622                              \\ 
\hline
         & 10                    & 153.0                & 423,873                            & infeasible                                  \\
3        & 50                    & 794.1                & 421,689                            & 436,608                              \\
         & 100                   & 1620                 & 408,388                            & 421,829                              \\ 
\hline
         & 10                    & 82.52                & 470,538                            & 470,538                              \\
4        & 50                    & 419.9                & 470,538                            & 470,538                              \\
         & 100                   & 842.3                & 470,538                            & 470,538                              \\ 
\hline
         & 10                    & 75.36                & 451,269                            & 470,869                              \\
5        & 50                    & 373.9                & 445,230                            & 458,624                              \\
         & 100                   & 748.6                & 445,230                            & 450,841                              \\ 
\hline
         & 10                    & 96.64                & 305,457                            & 337,105                              \\
6        & 50                    & 476.7                & 301,583                            & 314,667                              \\
         & 100                   & 949.7                & 298,414                            & 306,893                              \\ 
\hline
         & 10                    & 172.3                & 302,791                            & 358,768                              \\
7        & 50                    & 881.8                & 298,361                            & 313,493                              \\
         & 100                   & 1772                 & 298,361                            & 310,755                              \\ 
\hline
         & 10                    & 94.21                & 355,767                            & 395,748                              \\
8        & 50                    & 466.1                & 350,294                            & 363,609                              \\
         & 100                   & 934.3                & 348,585                            & 359,040                              \\
\hline\hline
\end{tabular}
\end{table}

\clearpage
\section*{Acknowledgments}
This material is based upon work supported by the National Science Foundation under Grant No. 1942559.

\bibliography{sample}
\end{document}